\documentclass{elsarticle}

\usepackage{amssymb,amsfonts,verbatim,graphicx,mathtools,
mathrsfs,commath,float,tikz-cd,amsthm,setspace}
\usetikzlibrary{positioning}
\usetikzlibrary{arrows}
\usepackage{caption}
\usepackage{subcaption}
\usepackage{tikz}
\usepackage{enumitem}
\usepackage[T1]{fontenc}
\usepackage{wrapfig}
\usepackage[hidelinks]{hyperref}
\usepackage{xfrac}

\usepackage[margin=1.2in]{geometry}

\newtheorem{ut}{Theorem}[section]
\newtheorem*{ut6}{Theorem 1.5}
\newtheorem{ul}[ut]{Proposition}
\newtheorem{uc}[ut]{Corollary}
\newtheorem*{uc7}{Corollary 1.6}

\newtheorem{uq}{Question} 
 
\theoremstyle{definition}

\theoremstyle{remark}
\newtheorem{ur}{Remark}

\makeatletter
\g@addto@macro\th@definition{\thm@headpunct{\textnormal{.}}}
\makeatother

\makeatletter
\def\ps@pprintTitle{%
 \let\@oddhead\@empty
 \let\@evenhead\@empty
 \def\@oddfoot{\centerline{\thepage}}%
 \let\@evenfoot\@oddfoot}
\makeatother

\DeclareMathOperator{\cnt}{cnt}
\DeclareMathOperator{\fcs}{fcs}

\DeclareMathOperator{\sng}{sng}

\numberwithin{equation}{section}
\numberwithin{ucl}{section}
\date{\today}
\makeatletter
\newcommand{\AlignFootnote}[1]{%
    \ifmeasuring@
    \else
        \footnote{#1}%
    \fi
}
\makeatother

\usetikzlibrary{calc}

\newcommand\doverline[1]{%
\tikz[baseline=(nodeAnchor.base)]{
    \node[inner sep=0] (nodeAnchor) {$#1$}; 
    \draw[line width=0.1ex,line cap=square] 
        ($(nodeAnchor.north west)+(0.0em,0.2ex)$) 
            --
        ($(nodeAnchor.north east)+(0.0em,0.2ex)$) 
        ($(nodeAnchor.north west)+(0.0em,0.5ex)$) 
            --
        ($(nodeAnchor.north east)+(0.0em,0.5ex)$) 
    ;
}}

\begin{document}

\begin{abstract}Suppose $Y$ is a continuum, $x\in Y$, and $X$ is the union of all nowhere dense subcontinua of $Y$ containing $x$.    Suppose further that there exists $y\in Y$ such that every connected subset of $X$ limiting to $y$ is dense in $X$. And, suppose  $X$ is dense in $Y$. We prove $X$ is homeomorphic to a composant of an indecomposable continuum, even though $Y$ may be decomposable. An example establishing the latter was given by Christopher Mouron and Norberto Ordoñez in 2016.   If $Y$ is chainable or, more generally, an inverse limit of identical topological graphs, then   we show $Y$ is  indecomposable and $X$ is a composant of $Y$. For homogeneous continua we explore similar problems which are related to a 2007 question of Janusz Prajs and Keith Whittington.




 \end{abstract}

\begin{frontmatter}

\title{Singularities of meager composants and filament composants}
\author{David Sumner Lipham}
\ead{dsl0003@auburn.edu}
\address{Department of Mathematics, Auburn University, Auburn, AL 36849}
\begin{keyword}continuum; meager composant; filament composant; singular; strongly indecomposable; homogeneous
\MSC[2010]
54F15, 54D35, 54H15
\end{keyword}

\end{frontmatter}

\section{Introduction}\label{s1}


\subsection{Terminology}
By a \textit{continuum} (plural form \textit{continua}) we shall mean a connected compact metrizable space with more than one point.   A continuum $Y$ is \textit{decomposable} if there are two proper subcontinua $H,K\subsetneq Y$ such that $Y=H\cup K$.   If $Y$ is not decomposable, then $Y$ is \textit{indecomposable}. We will say, more generally, that  a connected space $X$ is \textit{indecomposable} if $X$ cannot be written as the union of two proper closed connected subsets.  Equivalently, $X$ is indecomposable if $X$ is the only closed connected subset of $X$ with non-void interior \cite[\S48 V Theorem 2]{kur}.

A connected space $X$ is \textit{strongly indecomposable} if for every two non-empty disjoint open sets $U$ and $V$ there are two disjoint closed sets $A$ and $B$ such that $X\setminus U=A\cup B$, $A\cap V\neq\varnothing$, and $B\cap V\neq\varnothing$.  This term was introduced by the author in \cite{lip}. Strong indecomposability requires that the quasi-components of proper closed subsets of $X$ are nowhere dense in $X$, whereas indecomposability only requires that the connected components of proper closed subsets of $X$ are nowhere dense in $X$. 

Let $Y$ be a continuum and $x\in Y$.  The \textit{composant} of $x$ in $Y$ is the union of all proper subcontinua of $Y$ containing $x$.
Following  \cite{mou}, the   \textit{meager composant} of $x$ in $Y$ is the union of all nowhere dense subcontinua of $Y$ containing $x$.   A subcontinuum $K$ of $Y$ is said to be \textit{filament} if there exists a neighborhood of $K$ in which the connected component of $K$ is nowhere dense  \cite{pra}.  The   \textit{filament composant} of $x$ in $Y$ is the union of all filament subcontinua of $Y$ containing $x$.

 Given a  connected subset $X$ of a continuum $Y$,  and a point $y\in Y$, then $X$ is said to be \textit{singular with respect to $y$} if $\overline C=X$ for every connected $C\subseteq X$ with  $y\in \doverline{{C}}$.\footnote{Whenever $Y$ is a space of which $X$ is a subspace, and $A\subseteq X$, then we write $\overline A$ for the closure of $A$ in $X$, and $\doverline{A}$ for the closure of $A$ in $Y$.}  If there exists $y\in Y$ such that $X$ is singular with respect to $y$, then $X$ is \textit{singular} in $Y$.  And $X$ is \textit{singular dense} in $Y$ if  $X$ is both singular and dense in $Y$. This formulation is easily seen to be equivalent to the one in \cite{mou}.  

 
 A subset $X$ of a continuum $Y$ is called a \textit{filament set} if each continuum in $X$ is a filament subcontinuum of $Y$.  A continuum is \textit{filament additive} if  the union of every two intersecting filament subcontinua is filament \cite{pra2}. This property implies  the filament composants partition the continuum into pairwise disjoint sets. In homogeneous continua, filament additivity is equivalent to filament composants being filament sets  \cite[Corollary 3.6]{pra2}.\footnote{A continuum  is \textit{homogeneous} if for every two points $x$ and $y$ in the space there is a homeomorphism which maps the entire space onto itself, and maps $x$ to $y$.}   A continuum $Y$ is \textit{filamentable} if there is a filament subcontinuum $L\subseteq Y$ such that $Y\setminus L$ is a filament set.   

\subsection{Motivation and Summary of Results}There is no difference among composants,  meager composants, and filament composants  in an indecomposable continuum; \cite[Exercise 6.19]{nad} and \cite[Proposition 1.9]{pra}. And if $Y$ is an indecomposable continuum, and  $X$ is any composant of $Y$, then $X$ is singular with respect to each point of the dense $G_\delta$-set  $Y\setminus X$ \cite[Theorems 11.15 \& 11.17]{nad}.  The composant $X$ is also necessarily dense in $Y$ \cite[Theorem 5.4]{nad}.  

By contrast, an example in \cite[Section 5]{mou} shows the first two types of composants can differ quite dramatically inside of a decomposable continuum.   There was constructed a plane continuum $\mathfrak Y$ with only one traditional composant versus  uncountably many meager composants, each singular dense.

\begin{figure}[H]
  \centering
  \includegraphics[scale=1.3]{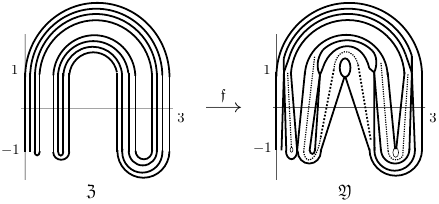} 
 \caption{The ``bucket-handle'' $\mathfrak Z$ and the Mouron-Ordoñez continuum $\mathfrak Y$ \textit{(graphics extracted from \cite[\href{ https://www.sciencedirect.com/science/article/pii/S0166864116301675}{Figure 3}]{mou})}.  Every meager composant of $\mathfrak Y$ is singular dense, but $\mathfrak Y$ is decomposable.}
 \label{f1}
  \end{figure}
  
That $\mathfrak Y$ has only one composant follows easily from the fact that $\mathfrak Y\cap ([0,3]\times [-1,1])$   is connected.\footnote{The continuum $\mathfrak Y\cap ([0,3]\times [-1,1])$ is  known as the ``Cajun accordion'' \cite{rog}.}  The (singular dense) meager composants of $\mathfrak Y$, on the other hand, are in one-to-one correspondence with the composants of the indecomposable ``bucket-handle'' continuum $\mathfrak Z$.   There is a continuous surjection $\mathfrak f:\mathfrak Z\to \mathfrak Y$ witnessing this fact. 
  Letting $\mathfrak X\subseteq \mathfrak Y$ be the image of the $\langle 0,-1\rangle$ endpoint composant of $\mathfrak Z$, we can see that $\mathfrak f\restriction (\mathfrak Z\setminus\mathfrak f^{-1}[\mathfrak X])$ is a homeomorphism (the sets $\mathfrak f^{-1} [\mathfrak X]\subseteq \mathfrak Z$ and $\mathfrak X\subseteq \mathfrak Y$ are indicated by the solid lines in Figure \ref{f1}). In particular,   every meager composant of $\mathfrak Y$ other than $\mathfrak X$ is homeomorphic to a composant of $\mathfrak Z$.\footnote{And these non-endpoint composants of $\mathfrak Z$ are mutually homeomorphic \cite{band}.}    By the first theorem of this paper, $\mathfrak X$ is also homeomorphic to a traditional composant.

\begin{ut}\label{t1}Every singular dense meager composant is homeomorphic to a composant of an indecomposable continuum.\end{ut}

Proving Theorem \ref{t1}  will  demonstrate that if  $Y$ is a continuum and $X$ is a singular dense meager composant of $Y$, then there is an indecomposable continuum $Z$ and homeomorphic embedding $\xi:X\hookrightarrow Z$ such that $Z$ has the same dimension as $X$; $\xi[X]$ is a composant of $Z$; and  there is a  mapping $f:Z\to Y$ such that $f\restriction \xi[X]=\xi^{-1}$ is a homeomorphism onto $X$.   This has  the following corollary.

\begin{uc}\label{cor3}If $Y$ is a continuum with  a singular dense meager composant, then each meager composant of $Y$ is dense.\end{uc}

Singularity is critical to Corollary \ref{cor3}.  For instance,  $[-1,0]^2\cup \mathfrak Z$ has both dense and non-dense meager composants. 

The next result applies to  chainable and circularly-chainable continua.  It builds on a  theorem of Mouron   \cite[Theorem 32]{mou2} stating that such a continuum is indecomposable if a sequence of disjoint subcontinua converges to the entire space in the Hausdorff metric.  

\begin{uc}Let $Y$ be a continuum which is the inverse limit of mutually homeomorphic topological graphs.  Then  $Y$ is indecomposable if and only if some (every) meager composant of $Y$ is singular dense.\label{t4}\end{uc}


Thus, there is no graph-like continuum like  $\mathfrak Y$.  Observe also that the meager composants of $\mathfrak Y$ are not filament sets.  For instance, the subcontinuum $\{0\}\times [-1,1]$ is non-filament. 

 \begin{ut}\label{18}A continuum $Y$  is indecomposable if and only if $Y$  has a meager composant which is also a singular dense filament set.\end{ut}

 The final results concern homogeneous continua.    
 
 \begin{ut}\label{t6}A homogeneous continuum $Y$ is indecomposable if and only if $Y$ is  filament additive, filamentable, and has singular dense filament composants.\end{ut}

\begin{uc}Let $Y$ be a filament additive, filamentable, homogeneous continuum with dense filament composants.  Then $Y$ is indecomposable if and only if the filament composants of $Y$ are indecomposable.\label{c7} \end{uc}

The sharpness of  the last two results  is evidenced by the product of a  circle with a (non-circle) solenoid. That continuum is  homogeneous, filament additive, filamentable, and decomposable. Its filament composants are products of the circle with composants of the solenoid \cite[Theorem 4.4]{pra2}.  These sets are dense, but are  neither singular nor indecomposable.

\section{Properties of indecomposable meager composants}

We begin by showing singular dense connecta are indecomposable.

\begin{ul}\label{p1}Let $X$ be a connected subset of a continuum $Y$. If $X$ is singular dense in $Y$, then $X$ is indecomposable.\end{ul}

\begin{proof}Suppose $X$ is singular dense in $Y$.  Let $y\in Y$ be such that $X$ is singular with respect to $y$.   For a contradiction suppose $X$ is  the union of two proper closed connected subsets $H$ and $K$.  By   $\doverline X=Y$ we  have  $y\in \doverline{{H}}$ or $y\in \doverline{{K}}$.  Neither ${\overline{H}}$ nor ${\overline{K}}$ is equal to $X$, so this contradicts $X$ being singular with respect to $y$. Therefore $X$ is indecomposable.\end{proof}

\begin{ur} \cite[Theorem 7.6]{mou} says that every irreducible continuum with a singular dense meager composant is indecomposable.  The proof uses the idea of  minimal decompositions.  Alternatively, combine Proposition \ref{p1} with \cite[Theorem 6(iii)]{lip}.  These results show that if $Y$ is an irreducible continuum in which any connected set is singular dense, then $Y$ is indecomposable. In particular, a continuum is indecomposable if (and only if) its composants are singular.
 \end{ur}

The next proposition shows that meager composants partition a continuum into pairwise disjoint sets (cf. \cite[Proposition 2.5]{mou}).

\begin{ul}\label{p2}Let $X$ be the meager composant of a point $x$ in a continuum $Y$.  Let $K$ be a nowhere dense subcontinuum of $Y$. If $K\cap X\neq\varnothing$, then $K\subseteq X$.\end{ul}

\begin{proof}Suppose $K\cap X\neq\varnothing$.   Let $x'\in K\cap X$.  There is a nowhere dense subcontinuum $L\supseteq \{x,x'\}$. Then $K\cup L$ is a nowhere dense subcontinuum of $Y$ containing $x$, so $K\subseteq X$. \end{proof}

For every topological space $A$ and point $x\in A$, we let $\cnt(x,A)$ denote the connected component of $x$ in $A$.  That is,  $\cnt(x,A)=\bigcup\{C\subseteq A:C\text{ is connected and }x\in C\}.$  When $A$ is a subset of a topological space $X$, then $A$ is always given the subspace topology.

\begin{ul}\label{p3}Let $X$ be a  meager composant of a continuum $Y$.   If $X$ is indecomposable, then:

\begin{enumerate}
\item[\textnormal{i.}]  every proper closed connected subset of $X$ is compact; 

\item[\textnormal{ii.}]   either $X$ is compact or $X$ is of the first category of Baire;  and

\item[\textnormal{iii.}] for every  $X$-closed  set $A\subseteq X$, the component decomposition $\mathfrak A:=\{\cnt(x,A):x\in X\}$ is metrizable and zero-dimensional.\footnote{The set $\mathfrak A$ is given  the quotient topology; $\mathfrak U$ is open (closed) in $\mathfrak A$ if and only if $\bigcup \mathfrak U$ is open (closed) in $A$. We say that a  space is \textit{zero-dimensional} if it has a basis of clopen (closed-and-open) sets.} 
\end{enumerate}\end{ul}

\begin{proof}Suppose  $X$ is indecomposable.

(i): Let $C$ be a proper closed connected subset of $X$.  Then   $\doverline{{C}}$ is a nowhere dense subcontinuum of $Y$ by indecomposability of $X$.  By Proposition \ref{p2}, $\doverline{{C}}\subseteq X$, so $\doverline{{C}}=\overline C=C$ is compact. 

(ii): Suppose $X$ is non-compact.  Let $x\in X$, and let $\{U_n:n<\omega\}$ be a basis for $X\setminus \{x\}$ consisting  of non-empty open sets.  Clearly $X\supseteq \bigcup\{\cnt(x ,X\setminus U_n):n<\omega\}$. Conversely, let $x'$ be any point in $X$.  There is a continuum $L\subseteq X$ with $\{x,x'\}\subseteq L$.  We know $L\neq X$ because $X$ is not compact, so there exists $n<\omega$ such that $L\cap U_n=\varnothing$. Then $x'\in \cnt(x,X\setminus U_n)$.  This shows   $X\subseteq \bigcup\{\cnt(x ,X\setminus U_n):n<\omega\}$.  Thus $X= \bigcup\{\cnt(x ,X\setminus U_n):n<\omega\}$. Each $\cnt(x, X\setminus U_n)$ is closed and nowhere dense by indecomposability of $X$. Therefore $X$ is of the first category of Baire.




(iii): Let $A$ be a closed subset of $X$. If $A=X$, then  $\mathfrak A=\{X\}$ is clearly metrizable and zero-dimensional.  Let us assume for the remainder of the proof that  $A\neq X$.    First we will show $\mathfrak A$ is metrizable. This will  be useful in proving $\mathfrak A$ is zero-dimensional.

Let $\varphi:A\to \mathfrak A$ be the canonical epimorphism defined by $\varphi(x)=\cnt(x,A)$.  To prove $\mathfrak A$  is metrizable, it suffices to show  $\varphi$  is perfect \cite[Theorem 4.2.13]{eng}.  Well,   each member of $\mathfrak A$ is compact by Proposition   \ref{p3}.i. It remains  to show $\varphi$ is closed.  To that end, let $C\subseteq A$ be closed.  To prove $\varphi[C]=\{A'\in \mathfrak A:A'\cap C\neq\varnothing\}$ is closed in $\mathfrak A$ we must show $\bigcup\varphi[C]$ is closed in $X$.

 Suppose $x\in \overline{\bigcup\varphi[C]}$. Then there exists $(x_n)\in (\bigcup\varphi[C])^\omega$ such that $x_n\to x$.  For each $n<\omega$ let $A_n=\cnt(x_n,A)$.

\textit{Case 1: A subsequence of $(A_n)$ has connected union.} Let  $(A_{n_k})$   be  a subsequence whose union is connected. Then  by maximality of the connected component $A_{n_0}$ we have $x_{n_k}\in A_{n_k}=A_{n_{0}}$ for each $k<\omega$. Since $A_{n_{0}}$ is closed we have  $x\in A_{n_{0}}\subseteq \bigcup\varphi[C]$.  

\textit{Case 2: No subsequence of $(A_n)$ has connected union.}  By     compactness of the hyperspace $K(Y)$ \cite[Corollary  4.13]{nad},  the sequence $(A_n)\in [K(Y)]^\omega$ has an accumulation point $M\in K(Y)$.\footnote{Here $K(Y)$ is the set of non-empty compact subsets of $Y$  equipped with the Vietoris topology \cite[2.7.20]{eng}.  For metrizable $Y$, the Vietoris topology  coincides with the topology generated by any Hausdorff metric \cite[4.5.23]{eng}.} 
  Necessarily, $x\in M\subseteq \doverline A$ and $M$ is a continuum \cite[Corollary 4.18]{nad}. $M$ is also nowhere dense in $Y$. For suppose otherwise that $M^\mathrm{o}\neq\varnothing$.  Then $\Omega\coloneqq \{n<\omega:A_n\cap M^\mathrm{o}\neq\varnothing\}$ is infinite. By hypothesis $\Sigma \coloneqq \bigcup \{A_n:n\in \Omega\}$ is not connected. So there are $Y$-open sets $U$ and $V$ such that $U\cap  \Sigma\neq\varnothing$, $V\cap \Sigma\neq\varnothing$, $U\cap V\cap \Sigma=\varnothing$, and   $\Sigma\subseteq U\cup V$.   The Vietoris open set $\{K\in K(Y):K\cap U\cap M^\mathrm{o}\neq\varnothing\text{ and }K\cap V\cap M^\mathrm{o}\neq\varnothing\}$ contains $M$ as an element but has empty intersection with $\{A_n:n<\omega\}$. This is a contradiction. Therefore $M^\mathrm{o}=\varnothing$, i.e. $M$ is nowhere dense.  By Proposition \ref{p2},   $M\subseteq A$. Further, $M\cap C\neq\varnothing$ (otherwise, letting $W$ be a $Y$-open set such that $W\cap X=X\setminus C$ we find that $M\subseteq W$ yet $W$ contains no $A_n$ because $A_n\cap C\neq\varnothing$).  Since $M\subseteq \cnt(x,A)$, this implies   $ \cnt(x,A)\in \varphi[C]$.  So $x\in \bigcup\varphi[C]$. 

In each of the two possible cases we found $x\in \bigcup\varphi[C]$.  Therefore $\bigcup\varphi[C]=\overline{\bigcup\varphi[C]}$, so that  $\varphi[C]$ is closed in $\mathfrak A$.  We conclude that $\varphi$ is perfect, so $\mathfrak A$ is metrizable. 

Next we show $\mathfrak A$ is zero-dimensional.  Note that if $X$ is compact then so is $A$, and in this case $\mathfrak A$ is  already known to be zero-dimensional \cite[Theorem 6.2.24]{eng}.   So assume  $X$ is non-compact.      
Let $U$ be an open subset of $Y$ such that $X\setminus U=A$.  By Proposition  \ref{p3}.ii there exists $y\in U\setminus X$.  Let $\varepsilon>0$ such that $B(y,\varepsilon)\subseteq U$, and for each $n<\omega$ put $W_n=B(y,\varepsilon/2^n)$. 

Fix   $x'\in A$. For each $x\in A$ there is a nowhere dense continuum $L\subseteq X$ such that $\{x,x'\}\subseteq L$. Since $y\notin L$ there exists $n<\omega$ such that $W_n\cap L=\varnothing$.  Then $x\in \cnt(x',X\setminus W_n)$. This shows
\begin{equation}\label{e21}A\subseteq \bigcup\{\cnt(x',X\setminus W_n):n<\omega\}.\end{equation} 

We also claim that \begin{equation}\label{e22} \cnt(x,A)=\cnt(x,K_n) \text{ for every }   x\in K_n\coloneqq  \cnt(x',X\setminus W_n)\cap A.\end{equation} Well, suppose $x\in K_n$. Then  $\cnt(x,A)\supseteq  \cnt(x,K_n)$ because $A\supseteq K_n$. Conversely,    $\cnt(x,A)\subseteq\cnt(x,X\setminus W_{n})= \cnt(x',X\setminus W_{n})$ implies   $$\cnt(x,A)=\cnt(x,\cnt(x,A)\cap A)\subseteq \cnt(x,\cnt(x',X\setminus W_{n})\cap A)=\cnt(x,K_n).$$ Therefore $\cnt(x,A)=\cnt(x,K_n)$.

For each $n<\omega$  let $\mathfrak K_n=\{\cnt(x,K_n):x\in K_n\}$ be the component decomposition of $K_n$. By (\ref{e21}) and (\ref{e22}) we have  $$\mathfrak A=\bigcup \{\mathfrak K_n:n<\omega\}.$$ 
Endow  the sets  $\mathfrak A$  and $\mathfrak K_n$  with the quotient topologies relative to $A$ and $K_n$, respectively,  and observe   that each $\mathfrak K_n$ is a subspace of  $\mathfrak A$.  For if $\mathfrak S$ is any subset of $\mathfrak K_n$, then:
\begin{align*}
\mathfrak S\text{ is closed in }\mathfrak K_n &\Leftrightarrow \bigcup \mathfrak S \text{ is closed in }K_n\\ 
&\Leftrightarrow \bigcup \mathfrak S\text{ is closed in }A \\
& \Leftrightarrow \mathfrak S\text{ is closed in }\mathfrak A.
\end{align*}
The first equivalence is the definition of the quotient topology on $\mathfrak K_n$.  The second equivalence holds because $K_n$ is a closed  subset of $A$.  The third holds by the inclusion $\mathfrak S\subseteq \mathfrak K_n\subseteq \mathfrak A$ and the definition of the quotient topology on $\mathfrak A$. 

Note that $K_n$  is compact by Proposition  \ref{p3}.i, so  $\mathfrak K_n$ is  zero-dimensional by  \cite[Theorem 6.2.24]{eng}.  Thus,      $\mathfrak A$ is a separable metrizable union of countably many closed (compact) zero-dimensional subspaces. By \cite[Theorem 1.3.1]{eng2},   $\mathfrak A$ is zero-dimensional. \end{proof}

The next proposition shows  that singularity of dense  meager composants can be expressed using various   familiar properties of  connected sets.  And for dense meager composants, being singular with respect to one point implies having a full complementary set of singularities.

   \begin{ul}\label{p5}Let  $X$ be a meager composant of a continuum $Y$.  If $X$ is dense in $Y$, then  the  following are equivalent: 
 \begin{enumerate}
\item[\textnormal{i.}] $X$ is indecomposable; 
  \item[\textnormal{ii.}] $X$ is strongly indecomposable; 
    \item[\textnormal{iii.}] there exists $y\in Y$ such that $X$ is singular  with respect to $y$;
       \item[\textnormal{iv.}] $Y\setminus X\neq\varnothing$ and $X$ is singular with respect to each point of $Y\setminus X$;
        \item[\textnormal{v.}] there exists $y\in Y$ such that the connected set $X\cup \{y\}$ is irreducible;
    \item[\textnormal{vi.}] $Y\setminus X\neq\varnothing$ and  $X\cup \{y\}$ is irreducible for every $y\in Y\setminus X$;
       \item[\textnormal{vii.}]  $X\cup \{y\}$ is indecomposable for every $y\in Y$.
  \end{enumerate}
  \end{ul}
  
\begin{proof}Suppose $\doverline X=Y$. 

First we will show  (ii)$\Rightarrow$(i)$\Rightarrow$(vi)$\Rightarrow$(v)$\Rightarrow$(iii)$\Rightarrow$(i)$\Rightarrow$(ii), establishing the equivalence of all items other than (iv) and (vii). Then, to incorporate (iv) and (vii)  we  will prove (vi)$\Rightarrow$(iv)$\Rightarrow$(iii) and (vi)$\Rightarrow$(vii)$\Rightarrow$(i).
  
(ii)$\Rightarrow$(i): Fairly obvious; see the second paragraph of \cite[Section 2]{lip}.
  
(i)$\Rightarrow$(vi): Suppose $X$ is indecomposable.  Then $X\neq Y$ because every indecomposable continuum has more than one  meager composant.  Let $y\in Y\setminus X$ and fix $x\in X$.  For a contradiction suppose $X\cup \{y\}$ is reducible between $x$ and $y$.  Let $C\supseteq \{x,y\}$ be a proper closed connected subset of $X\cup \{y\}$.   By Proposition \ref{p3}.iii the decomposition of $C\cap X$ into connected components is  metrizable and zero-dimensional, so there is a decreasing sequence of  $(C\cap X)$-clopen sets $E_0\supseteq E_1\supseteq E_2\supseteq...$ such that $\cnt(x,C\cap X)=\bigcap\{E_n:n<\omega\}$. Each  $E_n\cup \{y\}$ is connected, so  $K:=\bigcap\{\doverline{E_n\cup \{y\}}:n<\omega\}$ is the intersection of a decreasing sequence of continua.  Then $K$ is a continuum.    Further,  $K$ has non-empty interior because $X$ is a meager composant of $Y$, $\{x,y\}\subseteq K$, and $y\notin X$.  Let $U$ be a non-empty $Y$-open set such that $U\subseteq \doverline{E_n\cup \{y\}}$ for each $n<\omega$.  Then $U\cap X\subseteq \doverline{E_n\cup \{y\}}\cap X=E_n$ for each $n<\omega$. So  $\cnt(x,C\cap X)$, which is a proper closed connected subset of $X$, contains the non-empty $X$-open set $U\cap X$. This contradicts indecomposability of $X$. Therefore $X\cup \{y\}$ is irreducible (between $x$ and $y$).
    
(vi)$\Rightarrow$(v): Trivial. 

(v)$\Rightarrow$(iii): Suppose $y\in Y$ is such that $X\cup \{y\}$ is irreducible.   For every two points  $x$ and $x'$ in $X$  there is a continuum $L\subseteq X$ which contains $\{x,x'\}$ and  is nowhere dense in $Y$. Since $X$ is dense in $Y$, we know $L$ is also  nowhere dense in $X$. Therefore $X$ is reducible, so there exists  $x\in X$ such that $X\cup \{y\}$ is irreducible between $x$ and $y$.  Let $C$ be any connected subset of $X$ such that $y\in \doverline C$.  There is a nowhere dense (in $X$) continuum $L\subseteq X$ such that $L\cap C\neq\varnothing$ and $x\in L$.  Then $\overline C\cup L\cup \{y\}$ is a proper closed connected subset of $X\cup \{y\}$ containing $x$ and $y$.   By irreducibility  it must be that  $\overline C\cup L\cup \{y\}=X\cup \{y\}$, whence $\overline C=X$ and   $X$ is singular with respect to $y$.

  (iii)$\Rightarrow$(i): Proposition \ref{p1}.
  
(i)$\Rightarrow$(ii): Suppose $X$ is indecomposable.  Toward showing $X$ is strongly indecomposable, let $U$ and $V$ be non-empty disjoint open subsets of $X$.  We will exhibit a relatively clopen  subset of $X\setminus U$ which intersects $V$ but does not contain $V$. Well, by indecomposability of $X$ there are two connected  components $A_0\neq A_1$ of $X\setminus U$  such that $A_0\cap V\neq\varnothing$ and $A_1\cap V\neq\varnothing$.  The component decomposition of $X\setminus U$ is zero-dimensional by Proposition \ref{p3}.iii.   So in $X\setminus U$ there is a clopen set which contains $A_0$ and misses $A_1$.    

 (vi)$\Rightarrow$(iv): Similar to (v)$\Rightarrow$(iii). 
 
 (iv)$\Rightarrow$(iii): Trivial.

(vi)$\Rightarrow$(vii): For a contradiction suppose (vi) and the negation of (vii). Let $y\in Y$ be such that $X\cup \{y\}$ is decomposable.  Let $H$ and $K$ be proper closed connected subsets of $X\cup \{y\}$ such that $H\cup K=X\cup \{y\}$.  We have already established (vi)$\Rightarrow$(i), so $y\in (H\cap K)\setminus X$. By (vi), $X\cup \{y\}$ is irreducible.  Since $X$ is reducible, this means $X\cup \{y\}$ is irreducible between some $x\in X$ and $y$.  But for each $x\in X$ one of the sets $H$ or $K$ will show that $X\cup \{y\}$ is reducible between $x$ and $y$. This is a contradiction.

(vii)$\Rightarrow$(i):  Trivial.\end{proof}

\begin{ur}From the proof of (i)$\Rightarrow$(ii) we see that indecomposable meager composants are strongly indecomposable.    It remains an open problem to determine whether there is an indecomposable connected set which is  \textit{not} strongly indecomposable, but we suspect there is such an example.  Some variations of this problem appear in \cite[Section 5]{lip}.   \end{ur}

\begin{ur}Regarding (i)$\Rightarrow$(vii), in   \cite[Example 1]{lip} there was shown to be a locally compact indecomposable connected plane set whose one-point compactification  is decomposable. On the other hand, Mary Ellen Rudin \cite{rud} proved: If $X$ is any connected plane set and $Y$ is the plane closure of $X$,  then (i)$\Rightarrow$(vii).   We noticed that Rudin's proof could be dramatically simplified  if  every indecomposable connected plane set were known to be strongly indecomposable.  This suggests proving (i)$\Rightarrow$(ii) in general, or just for connected plane sets, could be  difficult.\end{ur}

\begin{ur}The implication (i)$\Rightarrow$(v) holds more generally when $X$ is any connected set and $Y$ is any compactification of $X$ \cite[Theorem 3]{lip2}. But (i)$\Rightarrow$(vi) is generally false by the example  in the previous remark.\end{ur}

\section{Proof of Theorem \ref{t1}}

\noindent Suppose $X$ is a singular dense meager composant of a continuum $Y$. We will construct an indecomposable continuum $Z$ with a composant homeomorphic to $X$.

By  strong indecomposability of $X$ (Proposition \ref{p5}.ii) there is a homeomorphic embedding $\iota:X \hookrightarrow [0,1]^\omega$ such that $I\coloneqq \overline{\iota[X]}$ is an indecomposable continuum \cite[Theorem 9]{lip}.  Let $\Gamma:Y\hookrightarrow  [0,1]^\omega$ be a homeomorphic embedding of $Y$, and put $\gamma\coloneqq \Gamma\restriction X$.   Let $\pi_n:[0,1]^\omega\to [0,1]$ be the $n$-th coordinate projection.  There is a homeomorphic embedding $\xi:X\hookrightarrow [0,1]^\omega$ such that all of the maps
\begin{align*}
\varphi_n\coloneqq \pi_n\circ \iota\circ \xi^{-1}:\;&\xi[X]\to [0,1]\text{; and}\\
\psi_n\coloneqq \pi_n\circ \gamma\circ \xi^{-1}:\;&\xi[X]\to [0,1]
\end{align*}
 continuously extend to  \begin{equation}Z:=\overline{\xi[X]}.  \label{e41}\end{equation}
 By   \cite[Exercise 1.7.C]{eng3},  $\xi$ can even be constructed to obtain $\dim(Z)=\dim(X)$.

 For each $n<\omega$ let  $\Phi_n:Z\to [0,1]$  and $\Psi_n:Z\to [0,1]$ be the continuous extensions of $\varphi_n$ and $\psi_n$, respectively.  Define $\Phi:Z\to [0,1]^\omega$ by $\pi_n\circ \Phi=\Phi_n$, and likewise for $\Psi$. Then $\Phi$ maps onto $I$ and  $\Psi$ maps onto $\Gamma[Y]$.  Since $\Phi\restriction \xi[X]=\iota\circ \xi^{-1}$ and $\Psi\restriction \xi[X]=\gamma\circ \xi^{-1}$ are homeomorphisms and $\xi[X]$ is dense in $Z$ (\ref{e41}),   we have \begin{align}\Phi^{-1}[ \iota[X]]&=\xi[X];\text{ and}\label{e42}\\
 \Psi^{-1}[ \gamma[X]]&=\xi[X].\label{e43}\end{align}
\vspace{-7mm}
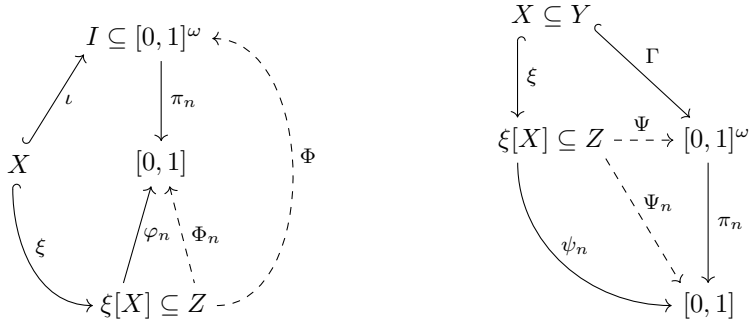
\begin{figure}[h]
\centering
     \begin{tikzpicture}[scale=1.11]
        \node (1) at (0.5,-1.5) {$X$};
        \node (2) at (2,0) {\textbf{$I\subseteq [0,1]^\omega$}};
        \node (3)  at (2.202,-1.5)  {\textbf{$[0,1]$}};
        \node (4)  at (1.8,-3.2)  {\hspace{7mm}$\xi[X]\subseteq Z$};
         \draw [right hook->] ([xshift=-.75mm,yshift=-.2mm]1.south) to [out=270,in=180]  node[right,above]{\hspace{2.5mm}\footnotesize$\xi$}([xshift=6.5mm,yshift=0mm]4.west) ;
         \draw [right hook->] ([xshift=1mm,yshift=0mm]1.north) --  node[right] {\hspace{0mm}\footnotesize$\iota$}([xshift=.9mm,yshift=.8mm]2.south west) ;
         \draw [->] ([xshift=2.02mm]2.south) --  node[right] {\hspace{-0mm}\footnotesize$\pi_n$}([xshift=0mm]3.north) ;
         \draw [dashed,<-] ([xshift=0mm]2.east) to [out=0,in=0]  node[right,above] {\hspace{4mm}\footnotesize$\Phi$}([xshift=0mm]4.east) ;
         \draw [dashed,<-] ([xshift=1mm]3.south) to  node[right] {\hspace{0mm}\footnotesize$\Phi_n$}([xshift=8mm]4.north) ;
         \draw[<-]  ([xshift=-1.2mm]3.south) -- node[right] {\hspace{-.4mm}\footnotesize$\varphi_n$}([xshift=-.5mm]4.north);
    \end{tikzpicture} \;\;\;\;\;\;\;\;\;\;\;\;\;\;\;\;\;\;\;\;
      \begin{tikzpicture}[scale=1.1]
        \node (1)  {$X\subseteq Y$};
        \node (2) at (0,-1.5) {\textbf{$\xi[X]\subseteq Z$}};
        \node (3)  at (2,-1.5)  {\textbf{$[0,1]^\omega$}};
        \node (4)  at (2,-3.5)  {\hspace{-2.04mm}\textbf{$[0,1]$}};
         \draw [right hook->] ([xshift=-1mm,yshift=1mm]1.south east) --  node[right,above]{\hspace{3mm}\footnotesize$\Gamma$}([xshift=2mm,yshift=0mm]3.north west) ;
         \draw [right hook->] ([xshift=1.9mm,yshift=0mm]1.south west) --  node[right] {\hspace{0mm}\footnotesize$\xi$}([xshift=-4.09mm,yshift=0mm]2.north) ;
         \draw [dashed,->] (2.east) --  node[left,above] {\hspace{-.5mm}\footnotesize$\Psi$}([xshift=0mm]3.west) ;
         \draw [dashed,->] ([yshift=.5mm,xshift=-1mm]2.south east) --  node[right,above] {\hspace{4mm}\footnotesize$\Psi_n$}([xshift=-1mm]4.north west) ;
         \draw [->] ([xshift=3.5mm]2.south west) to [out=270,in=180] node[right,above] {\hspace{3mm}\footnotesize$\psi_n$}([xshift=-2.06mm]4.west) ;
         \draw[->]  ([xshift=-1.06mm]3.south) -- node[right] {\footnotesize$\pi_n$}([xshift=-1.06mm]4.north);
    \end{tikzpicture}
\caption{Commutative diagrams for $\Phi$ and $\Psi$.}
\end{figure}

\textit{$Z$ is indecomposable}: By (\ref{e41}) and (\ref{e42}), $\Phi$ maps onto $I$ and maps proper subcontinua of $Z$ to proper subcontinua of $I$.  Indecomposability of $I$  therefore implies $Z$  is indecomposable. For if $Z$ were the union of two proper subcontinua $H$ and $K$, then $I$ would be the union of proper subcontinua $\Phi[H]$ and $\Phi[K]$.

\textit{$\xi[X]$ is  contained in a composant of $Z$}: 
 Since $X$ is singular dense in $Y$ we know that $X$ is not compact. Therefore $\xi[X]\neq Z$. Also $X$ is continuum-wise connected,  thus $\xi[X]$ is  contained in a composant of $Z$.

\textit{$\xi [X]$ contains  a   composant of $Z$}: Let $x\in X$, and let $N\ni \xi(x)$ be a proper subcontinuum of $Z$. We will show $N\subseteq \xi[X]$.  Well, since  $Z$ is indecomposable $N$ is nowhere dense in $Z$. Therefore $\Psi[N]$ is a nowhere dense subcontinuum of $\Gamma[Y]$. For otherwise,  $\Psi[N]$ contains a  $\gamma[X]$-open set $U\neq\varnothing$.  By  (\ref{e43}),  $N$ contains the non-empty $\xi[X]$-open set $\xi\circ\gamma^{-1}[U]$.  Since $N$ is closed in $Z$ and (\ref{e41}) holds,   this implies $N$ has non-void interior in $Z$, a contradiction.  Thus $\Psi[N]$ is nowhere dense.  Since $\gamma[X]$ is a meager composant of $\Gamma[Y]$  we have $\Psi[N]\subseteq \gamma[X]$.  By (\ref{e43}), $N\subseteq \xi[X]$.

The composants of $Z$ are pairwise disjoint, so the two containments show that $\xi[X]$ is equal to a composant of $Z$. This completes the proof of Theorem \ref{t1}. \qed

\begin{ur}Define $f= \Gamma^{-1}\circ \Psi$ for a surjection  $f:Z\to Y$  such that $f[\xi[X]]=X$. By monotone-light factorization \cite[Theorem 13.3]{nad},  $f$ is equal to a monotone mapping of $Z$ onto some continuum $M$, followed by a surjective mapping $l:M\to Y$ such that $l^{-1}\{y\}$ is totally disconnected for every $y\in Y$. We see that $M$ is also an indecomposable continuum containing $l^{-1}[X]\simeq X$ as a composant.  \end{ur}

\section{Proof of Corollary \ref{cor3}}

\noindent If continuum $Y$ has a singular dense meager composant, we have  shown that a continuum $Z$  maps onto $Y$ so that each meager composant of $Y$ contains the image of a composant of $Z$.  Each composant of $Z$ is dense, therefore each meager composant of $Y$ is dense. \qed

\begin{ur}We now see that if  $Y$ is a continuum with a singular dense meager composant, then $Y$ has at least two (disjoint) dense meager composants.  In particular, for each $x\in Y$  there exists $y\in Y$ such that the union of all continua in $Y\setminus \{x\}$ containing $y$  is dense in $Y$.   This is enough to imply $Y$ is indecomposable if  $Y$ is chainable \cite[Corollary 3.6]{oo}.     \end{ur}

\section{Proof of Corollary \ref{t4}}

\noindent We prove only the non-standard implication.  

Suppose $X$ is a singular dense meager composant of the graph-like continuum $Y$.   To show $Y$ is indecomposable, by \cite[Theorem 32]{mou2}  it suffices to show there is a sequence $(X_i)$ of pairwise disjoint continua in $Y$ such that $d_H(X_i,Y)\to 0$, where $d_H$ is the Hausdorff distance induced by a metric $d$ on $Y$.   It is unknown whether $Y$ must have infinitely many meager composants (see Question 1 in Section 7), but in any case the continua  $X_i$ can be selected from $X$.  

Identify $X$ with a composant of an indecomposable continuum $Z$ which maps onto $Y$   (Theorem \ref{t1}).  Let $\varrho_H$ be the Hausdorff metric generated by a metric $\varrho$ on $Z$, and define $\varrho_{\inf}(z,A)=\inf\{\varrho(x,z):z\in A\}$ for each $x\in X$ and $A\subseteq Z$.  Let $z\in Z\setminus X$, and let  $(x_n)\in X^\omega$ such that $x_n\to z$.   Recursively define $X_i$  as follows. Put  $X_0=\{x_0\}$. There exists a positive integer $n_1$ such that $\varrho(x_0,x_{n_1})>1/n_1$ and $\varrho_H(\cnt(x_{n_1},Z\setminus B_{\varrho}(x_0,1/n_1)),Z)<1.$ If there were no such integer, then  `boundary bumping' \cite[Lemma 6.1.25]{eng} and compactness of the hyperspace $K(Z)$ would reveal a proper subcontinuum of $Z$ containing both $x_0$ and $z$.  Set  $X_1=\cnt(x_{n_1},Z\setminus B_{\varrho}(x_0,1/n_1))$. 

Suppose $i>1$ and $X_j\subseteq X$ has been defined for each $j<i$.  No proper subcontinuum  of $Z$ containing $z$ also meets the compact set $X_0\cup X_1\cup ...\cup X_{i-1}$, so there is a sufficiently large integer $n_i$ such that $\varrho_{\inf}(x_{n_i},X_0\cup X_1\cup ...\cup X_{i-1})>1/n_i$ and $$\varrho_H(\cnt(x_{n_i},Z\setminus B_{\varrho_{\inf}}(X_0\cup X_1\cup ...\cup X_{i-1},1/n_i)),Z)<1/i.$$ Let $X_i=\cnt(x_{n_i},Z\setminus B_{\varrho_{\inf}}(X_0\cup X_1\cup ...\cup X_{i-1},1/n_i))$. 

The terms of the sequence $(X_i)$ are pairwise disjoint continua in $X$, and $\varrho_H(X_i,Z)\to 0$.  Since $Z$ maps continuously onto $Y$, we have $d_H(X_i,Y)\to 0$.  \qed


\section{Proof of Theorem \ref{18}}

\noindent Let $X$ be a meager composant of a continuum $Y$.

If $Y$ is indecomposable, then by elementary continuum theory $X$ is both a composant of $Y$ and a singular dense filament  subset of $Y$. 
 
 Now suppose $X$ is a singular dense filament set. Let $x\in X$, and let $A\subseteq Y$ be a minimal non-filament subcontinuum  containing $x$ provided by \cite[Corollary 1.13]{pra}.   Since $X$ is a filament set there exists $y\in  A\setminus X$.   Let $C$ be the  composant of $x$ in $A$. 
Then $C\subseteq X$ by minimality of $A$ and the fact that filament subcontinua are nowhere dense.  Further $\doverline C=A$ \cite[Exercise 5.20]{nad}, so $y\in \doverline C$. By Proposition \ref{p5},  $X$ is singular with respect to $y$, so $\overline C= X$. By density of $X$ in $Y$ we have $\doverline{{C}}= Y$. Therefore  $A=Y$, so  $Y$ is indecomposable. \qed

\section{Questions}
\begin{uq}\label{q1}Let $Y$ be a continuum with a singular dense meager composant. Must $Y$ have at least three meager composants?\end{uq}

Compare Question \ref{q1} with \cite[Problem 8.8]{mou} on whether there is a continuum $Y$ with a point $x$ such that the  meager composants of $Y$ are $\{x\}$ and $Y\setminus \{x\}$. A counterexample to Question \ref{q1}  would also have exactly two meager composants;  a dense first category $F_\sigma$-set and its complement, a dense $G_\delta$-set.  This $G_\delta$ would not be  $F_\sigma$, contrary to   \cite[Conjecture 8.4]{mou}. Therefore, we conjecture a positive answer to Question \ref{q1}. 

 We would  like to know if there is a homogeneous example like $\mathfrak Y$.

 \begin{uq}\label{q3}Is there a decomposable homogeneous  continuum with singular dense meager composants?\end{uq}

Analogous questions for filament composants are also of interest. 

 \begin{uq}\label{q2}Is there a decomposable continuum with singular dense filament composants? \end{uq}
 
\begin{uq}\label{q3}Is every filament additive homogeneous continuum  with singular dense filament composants  necessarily indecomposable?\end{uq}

The next section contains some results toward answering Question 4 in the affirmative.

\section{Filament singularities in filament additive homogeneous continua}

Here we will prove Theorem  \ref{t6}  and Corollary \ref{c7}, and show that Question \ref{q3}  is related to a question of Prajs \& Whittington.    Throughout this section, \textit{$Y$ is assumed to be a filament additive homogeneous continuum with dense filament composants.}

For each point $x\in Y$  let $\fcs(x)$ denote  the filament composant of $x$. Let $$\sng(x)=\{y\in Y:(\forall\text{ connected }C\subseteq \fcs(x)\text{ with }x\in C)(y\in \doverline{C}\Rightarrow \doverline C=Y)\}$$ be the set of filament singularities of $x$.   

\begin{ul}$y\in \sng(x)$ if and only if $\fcs(x)$ is singular with respect to $y$. \label{81} \end{ul}

\begin{proof}Suppose $\fcs(x)$ is singular with respect to $y$.  This means if $C$ is any connected subset of $\fcs(x)$  (with or without the base point $x$) and $y\in \doverline C$, then $\overline C=\fcs(x)$.  Further,  $\doverline C=Y$ by the standing assumption that $\fcs(x)$ is dense in $Y$. This shows $y\in \sng(x)$. Now suppose $\fcs(x)$ is \textit{not} singular with respect to $y$.  Then there is a connected set $C\subseteq \fcs(x)$ such that $y\in \doverline C$ and $\overline C\neq\fcs(x)$.   Let $x'\in C$, and let $L\supseteq \{x,x'\}$ be a filament subcontinuum of $Y$.  Then $C\cup L$ is a connected subset of $\fcs(x)$, $x\in C\cup L$, and $y\in \doverline{C\cup L}\neq Y$.  Thus $y\notin \sng(x)$.\end{proof}

Similar to Proposition \ref{p5}:

\begin{ul}\label{pi} The following are equivalent:

\begin{enumerate}
\item[\textnormal{i.}] $\sng(x)\neq\varnothing$;
\item[\textnormal{ii.}] $\fcs(x)$ is singular (as defined in Section \ref{s1});
\item[\textnormal{iii.}] $\fcs(x)$ is indecomposable.
\end{enumerate}\end{ul}

\begin{proof}(i)$\Rightarrow$(ii): Proposition \ref{81}. (ii)$\Rightarrow$(iii): Proposition \ref{p1}. (iii)$\Rightarrow$(i): Suppose $\fcs(x)$ is indecomposable. Let $\{V_n:n<\omega\}$ be a basis of non-empty open sets for $Y\setminus \{x\}$. By $\doverline{\fcs(x)}=Y$ and indecomposability of $X$ each $\cnt(x ,\fcs(x)\setminus V_n)$ is nowhere dense, so there exists  $$y\in Y\setminus \bigcup \big\{\doverline{{\cnt(x ,\fcs(x)\setminus V_n)}}:n<\omega\big\}.$$ 
Then  $y\in \sng(x)$.   \end{proof}

\begin{ul}If $x'\in \fcs(x)$, then $\sng(x)=\sng(x')$.\label{p8}\end{ul}

\begin{proof}Let $x'\in \fcs(x)$.  By filament additivity $\fcs(x)=\fcs(x')$.  So for every $y\in  Y$, $\fcs(x)$ is singular with respect to $y$ if and only if $\fcs(x')$ is  singular with respect to $y$.  By Proposition \ref{81} $\sng(x)=\sng(x')$. \end{proof}

A subcontinuum $A\subseteq Y$ is \textit{ample} if $\cnt(A,U)$ is a neighborhood of $A$ for each open set $U\supseteq A$ \cite{pra}.  In homogeneous continua, \textit{ample}  and  \textit{non-filament} are equivalent \cite[Proposition 2.3]{pra}.  Minimal ample subcontinua of $Y$ exist by \cite[Corollary 2.5]{pra}.

\begin{ul}\label{p66}Let $A$ be a minimal ample subcontinuum of $Y$.  If $A\neq Y$, then $A\cap \sng(x)=\varnothing$   for each $x\in A$.
\end{ul}

\begin{proof}Let $x\in A$. Let $C$ be the  composant of $x$ in $A$.  Then $C$ is connected, $C\subseteq \fcs(x)$,  and $\doverline{{C}}=A$.   Thus $A\neq Y$ implies  $A\cap \sng(x)=\varnothing$.  \end{proof}

\begin{ul}\label{t3}$Y$ is indecomposable if and only if $\fcs(x)\cup \sng(x)=Y$ for some (every) $x\in Y$.\end{ul}

\begin{proof}If $Y$ is indecomposable then the equation holds for each $x$ because the filament composants and traditional composants of $Y$ coincide and partition $Y$.  Conversely, if $Y$ is decomposable then there is a  minimal ample subcontinuum $A\neq  Y$.  Let $x\in A$.   Since $\fcs(x)$ is a filament set \cite[Corollary 3.6]{pra2}, $A\setminus \fcs(x)\neq\varnothing$. Also, $A\cap \sng(x)=\varnothing$ by Proposition \ref{p66}. Therefore  $A\setminus [\fcs(x)\cup \sng(x)]\neq\varnothing$, whence $\fcs(x)\cup \sng(x)\neq Y$.  Since $Y$ is homogeneous we have $\fcs(x)\cup \sng(x)\neq Y$ for every $x\in Y$.\end{proof}

By Propositions \ref{81} and \ref{t3}, $Y$ is indecomposable if and only if the filament composants of $Y$ are singular with respect to all points in their respective complements.  Density of filament composants is critical to this result.  For example, the circle of pseudoarcs is homogeneous, filament additive, and decomposable.  Its  filament composants are singular with respect to all points in their complements, but are not dense.

Let us now examine the non-singularity relation $$\neg\sng=\{\langle x,y\rangle \in Y^2:y\notin \sng(x)\}.$$  Write $\neg\sng\langle x,y\rangle$ for $\langle x,y\rangle\in \neg\sng$.

\begin{ul}$\neg\sng$ is an equivalence relation.\label{p7}\end{ul}

\begin{proof}We need to show $\neg\sng$ is reflexive, symmetric, and transitive.

\textit{Reflexive:}  The standing assumption $\doverline{\fcs(x)}=Y$ implies  $\{x\}\subsetneq \fcs(x)$. Therefore $\neg\sng\langle x,x\rangle$.

\textit{Symmetric:} Suppose $\neg\sng\langle x,y\rangle$.  We will show $\neg\sng\langle y,x\rangle$.  To that end, let $C\ni x$ be a connected subset of $\fcs(x)$ such that $y\in \doverline{{C}}\neq Y$.  Let $p\in Y\setminus \doverline{C}$, and let $\varepsilon>0$ such that $B(p,2\varepsilon)\cap \doverline{{C}}=\varnothing$. For each $n<\omega$: let $\delta_n$ be an Effros number\footnote{If $Y$ is a homogeneous continuum, then for every $\varepsilon>0$ there is a positive number $\delta$, called an \textit{Effros number} for $\varepsilon$, such that for each pair of points $x$ and $y$ with $d(x, y)<\delta$ there is an onto homeomorphism $h:Y\to Y$ such that $h(x)=y$ and $d(z, h(z)) <\varepsilon$ for each $z \in Y$. This is called the Effros Theorem. It follows from the more general \cite[Theorem 2]{eff}.} for $\varepsilon/2^n$;  let  $x_n\in C$ such that $d(x_n,y)<\delta_n$; and  let $h_n:Y\to Y$ be a surjective homeomorphism such that $h_n(x_n)=y$ and $d(z,h_n(z))<\varepsilon/2^n$ for all $z\in Y$. The connected set $E:=\bigcup \{h_n[C]:n<\omega\}$ shows  $\neg\sng\langle y,x\rangle$.   Indeed,  $y\in E$, and $E\subseteq \fcs(y)$ by filament additivity and the fact that homeomorphisms respect filament composants.  Further,  $x\in \doverline{{E}}$ because $h_n(x)\to x$, and $\doverline{{E}}\cap B(p,\varepsilon)=\varnothing$. Therefore $\neg\sng\langle y,x\rangle$.

\textit{Transitive:} Suppose $\neg\sng\langle x,y\rangle$ and $\neg\sng\langle y,z\rangle$. We will show $\neg\sng\langle x,z\rangle$.

If $\sng(x)=\varnothing$ then clearly $\neg\sng\langle x,z\rangle$.  Now suppose   $\sng(x)\neq\varnothing$.  Then  $\fcs(x)$ is indecomposable by Proposition \ref{pi}. By $\neg\sng\langle x,y\rangle$ there is a connected set $C\subseteq \fcs(x)$ such that $x\in C$ and $y\in \doverline{{C}}\neq Y$.  By $\neg\sng\langle y,z\rangle$ and symmetry of $\neg\sng$,  there is also a connected set $D\subseteq \fcs(z)$ such that $z\in D$ and  $y\in \doverline{{D}}\neq Y$.  Indecomposability of $\fcs(x)$ implies $\doverline C$ is nowhere dense, so $\doverline{C}\cup \doverline{D}\neq Y$.  There exists $q\in Y\setminus \doverline{C\cup D}$ and $\varepsilon>0$ such that $B(q,2\varepsilon)\cap \doverline{C\cup D}=\varnothing$.  For each $n<\omega$ let $\delta_n$ be an Effros number  for $\varepsilon/2^n$. Let $x_n\in C \cap B(y,\delta_n/2)$ and $z_n\in D\cap B(y,\delta_n/2)$,  so that $d(x_n,z_n)<\delta_n$. Let $h_n:Y\to Y$ be a surjective homeomorphism such that $h_n(z_n)=x_n$ and $d(w,h_n(w))<\varepsilon/2^n$ for all $w\in Y$.  The connected set $C\cup\bigcup\{ h_n[D]:n<\omega\}\subseteq \fcs(x)$  witnesses $\neg\sng\langle x,z\rangle$.  \end{proof}
   \begin{figure}[H]
  \centering
  \includegraphics[scale=.44]{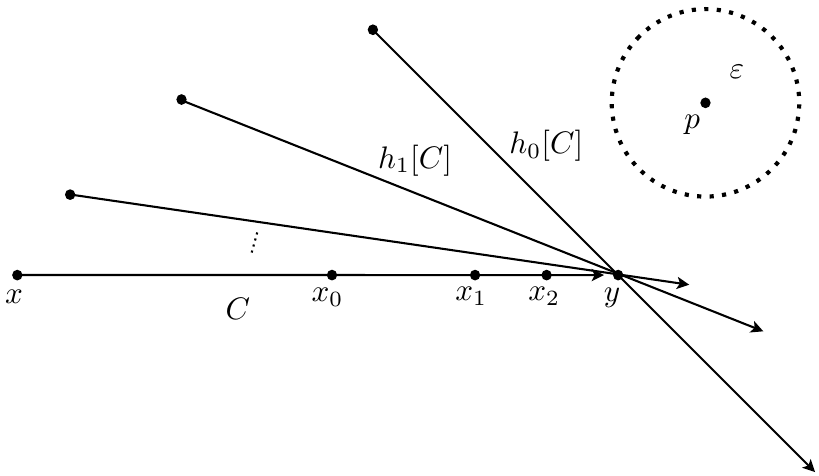}
  \hspace{5mm}
 \includegraphics[scale=.44]{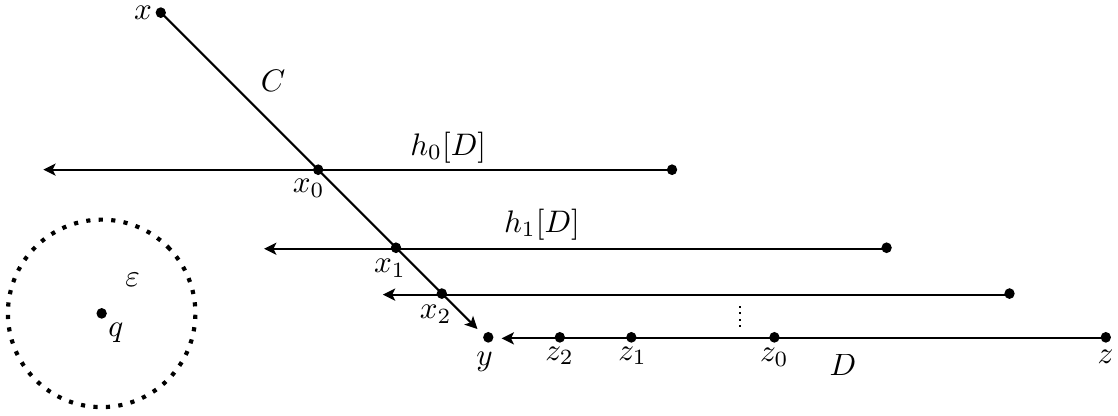}  
 \caption{Symmetry (left) and Transitivity (right) of $\neg\sng$ (Proposition \ref{p7}).}
 \label{f3}
  \end{figure}


\begin{ur}By Proposition  \ref{p7},  the non-singularity relation partitions $Y$ into pairwise disjoint sets.  By Proposition \ref{p8} and symmetry of $\neg\sng$, the partition $Y/\neg\sng=\{Y\setminus \sng(x) :x\in Y\}$ is  coarser than the partition of $Y$ into  filament composants. Likewise, each $\sng(x)$ is a union of filament composants.\end{ur}

\begin{ur}Proposition \ref{t3} says  $Y$ is indecomposable if and only if $\neg\sng=\fcs$, i.e. $$Y/\neg\sng=\{\fcs(x) :x\in Y\}.$$\end{ur}

\begin{ur}Question \ref{q3} asks whether $|Y/\neg\sng|>1$ implies $Y$ is indecomposable.\end{ur}


Let us now restate and prove the last two items from Section 1.2.

\begin{ut6}A homogeneous continuum $Y$ is indecomposable if and only if $Y$ is  filament additive, filamentable, and has singular dense filament composants.\end{ut6}

\begin{proof}Suppose the homogeneous continuum $Y$ is filament additive, filamentable, and has singular dense filament composants. Let $L$ be a filament subcontinuum of $Y$ such that $Y\setminus L$ is a filament set.    By singularity and  Proposition \ref{p7}, there exists exists $x\in Y$ such that $L\cap \neg\sng[x]=\varnothing$. Then $\neg\sng[x]=\bigcup\{\doverline C:C\subseteq \fcs(x)\text{ is connected, } x\in C \text{, and }\doverline C\neq Y\}$ is a continuum-wise connected filament set, so $\neg\sng[x]=\fcs(x)$. By Proposition \ref{t3}, $Y$ is indecomposable.\end{proof}

\begin{uc7}Let $Y$ be a filament additive, filamentable, homogeneous continuum with dense filament composants.  Then $Y$ is indecomposable if and only if the filament composants of $Y$ are indecomposable.\end{uc7}

\begin{proof}If the filament composants of $Y$ are indecomposable, then so is $Y$ by Proposition \ref{pi} and Theorem \ref{t6}. The converse follows from the fact that dense connected subsets of indecomposable continua are indecomposable. \end{proof}

\begin{ur}Based on Proposition \ref{p66}, the existence of one minimal ample subcontinuum of $Y$ which meets two non-singularity classes would imply $Y$ is indecomposable. Since each non-singularity class is a union of filament composants, a positive answer to the  question of Prajs and Whittington below would imply a positive answer to Question 4.  \end{ur}

\begin{uq}[{Question 7 in }\cite{pra2}] Let $Y$ be a filament additive homogeneous continuum with
dense filament composants. Does each ample subcontinuum of $Y$ intersect every
filament composant of $Y$?\end{uq}









\end{document}